\documentclass{amsart}

\usepackage{graphicx}
\usepackage{amsmath,mathrsfs}
\usepackage{amssymb,mathtools}
\usepackage{paralist}
\usepackage[latin1]{inputenc} 

\def\squareforqed{\hbox{\rlap{$\sqcap$}$\sqcup$}}
\def\qed{\ifmmode\squareforqed\else{\unskip\nobreak\hfil
\penalty50\hskip1em\null\nobreak\hfil\squareforqed
\parfillskip=0pt\finalhyphendemerits=0\endgraf}\fi}

\newcommand{\prt}{\partial}
\newcommand{\sub}{\subset}

\newcommand{\ria}{\rightarrow}

\newcommand{\lan}{\langle}
\newcommand{\ran}{\rangle}

\renewcommand{\phi}{\varphi}
\newcommand{\al}{\alpha}
\newcommand{\om}{\omega}
\newcommand{\Om}{\Omega}

\newcommand{\la}{\lambda}
\newcommand{\De}{\Delta}
\newcommand{\de}{\delta}
\newcommand{\ep}{\varepsilon}

\renewcommand{\th}{\theta}

\newcommand{\be}{\beta}
\newcommand{\ga}{\gamma}

\newcommand{\sig}{\sigma}

\newcommand{\mR}{{\mathbb{R}}}

\newcommand{\scp}[2]{{\big\lan {#1}\, , \, {#2}\big\ran}}

\newcommand{\defset}[1]{\left\{ #1 \right\}}

\newcommand{\abs}[1]{\left| #1 \right|}

\newcommand{\map}[5]{\begin{aligned}[b]#1: &&#2&\ria #3\\&&#4&\mapsto #5\end{aligned}}

\newcommand{\ignore}[1]{}

\newcommand{\Diff}{\mathrm{Diff}}

\newcommand{\id}{{\mathrm{id}}}

\newcommand{\bfone}{{\mathbf 1}}
\newcommand{\acos}{{\mathrm{arccos}}}

\newcommand{\sign}{{\mathrm{sign}}}

\newtheorem{theorem}{Theorem}

\graphicspath{{./Fig/}}

\title{Elastic distance between  curves under the metamorphosis viewpoint}
\author{Laurent Younes}

\begin{document}

\begin{abstract}
We provide a new angle and obtain new results on a class of metrics on length-normalized curves in $d$ dimensions, represented by their unit tangents  expressed as a function of arc-length, which are  functions from the unit interval to the $(d-1)$-dimensional unit sphere. These metrics are derived from the combined action of diffeomorphisms (change of parameters) and arc-length-dependent rotation acting on the tangent. Minimizing a Riemannian metric balancing a right-invariant metric on diffeomorphisms and an $L^2$ norm on the motion of tangents leads to a special case of ``metamorphosis'', which provides a general framework adapted to similar situations when Lie groups acts on Riemannian manifolds. Within this framework and using a  Sobolev norm with order 1 on the diffeomorphism group, we generalize previous results from the literature that provide explicit geodesic distances on parametrized curves. 

\end{abstract}

\maketitle

\section{Introduction}

There has been, over the past twenty years, a sizable amount of  work exploring elastic distances between plane curves and their computation using a square-root transformation mapping the space of curves into some standard infinite-dimensional Riemannian manifold.  In \cite{azencott1996distance,you96,you98} a  distance between parametrized plane curves was introduced, in which a transformation of the pair $(\phi, \theta)$ (involving a square root) placed the metric in a Hilbert space context, where $\phi$ is the parametrization and $\theta$ is the tangent angle (as functions of the arc-length). This distance can then be optimized with respect to $\phi$ to yield a metric between curves modulo reparametrization (a.k.a. unparametrized curves). 
Existence results for minimizers were then provided in \cite{ty00,ty01}.
Further analysis were made in the smooth case, with isometries with Stiefel and Grassman manifolds for closed curves and closed curves modulo rotations \cite{ymsm07}.  A different (but similar) approach was introduced in \cite{mio2007shape,ksmj04} and further developed in numerous papers, among which \cite{jksj07,kurtek2012elastic,samir2012gradient,su2014statistical,xie2013parallel,srivastava2016functional}, to provide a metric between curves, also using a square root to reduce to a Hilbert case.  More recently, the authors in \cite{bauer2014constructing} designed a different isometry applicable to a family of distances that includes the previous two.

In this paper, we reinterpret this line of work under the viewpoint of metamorphosis, which is described and developed in \cite{my01,ty05,ty05b,hty09,richardson2013computing,Richardson2016}. This reformulation will allow us to generalize previous results on the subject by placing them in a unified context. 

What we mean by elastic distances between curves are Riemannian metrics in spaces of parametrized curves that is, when evaluated at a smooth vector field along a curve, equivalent to the square norm of the derivative of this vector field with respect to the arc length. This is a small part of the range of Riemannian metrics that were considered in the literature. We refer to \cite{mm06a,mm07} for an extensive catalog and properties and to \cite{bauer2014constructing,bauer2017numerical} for more recent developments.

\section{Comparing Curves}
\index{Metamorphosis!plane curves}
\subsection{Metamorphosis on unit tangents}
\label{sec:meta.curves.1}

Metamorphosis describes  a general approach to build new Riemannian metrics on Riemannian manifolds acted upon by Lie groups \cite{my01,ty05,ty05b,hty09}. In a nutshell, letting $G$  be a Lie group acting on a manifold $M$, we define the action map $\pi: G\times M \to M$ by $\pi(\phi, \al) = \phi\cdot\al$ and the right action of $G$ on $G\times M$ by $(\phi, \al)\cdot \psi = (\phi\psi, \psi^{-1}\cdot \al)$. A metamorphosis is just a curve $t\mapsto (\phi(t), \al(t))$ in $G\times M$ and its image is $a(t) = \phi(t)\cdot \al(t)$.

Any right invariant metric on $G\times M$ specifies a unique metric on $M$ such that $\pi$ is a Riemannian submersion, and optimal metamorphoses associated with this metric are horizontal geodesics in $G\times M$. More precisely, if $R_\phi$ denotes the right translation by $\phi$, then a right-invariant metric on $G\times M$ must satisfy (letting $\id$ denote the identity element of $G$)
\[
\|(v, \xi)\|_{(\phi,\al)} = \|dR_{\phi^{-1}}(\phi,\al)\, (v, \xi)\|_{(\id, \phi\cdot \al)}
\]
and optimal metamorphoses are curves in $G\times M$ that minimize
\begin{equation}
\label{eq:meta.gen}
\int_0^1 \|dR_{\phi^{-1}(t)}(\phi(t),\al(t))\, (\prt_t \phi, \prt_t \al)\|^2_{(\id, \phi(t)\cdot \al(t))}\, dt
\end{equation}
with fixed initial conditions $\phi(0) = \id$,  $\al(0) = a_0$, $\phi(1)\cdot \al(1) = a_1$, $a_0, a_1\in M$. (The assumption that $\phi(0) = \id$ is without loss of generality because of right invariance.)\\

We  apply this principle in order to compare curves based on the orientation of their tangent. If $m: [0, L]\to \mR^d$  is parametrized with arc length, its normalized tangent is
$$\map{T^m}{[0,1]}{\mR^d}{s}{\prt_s m(Ls)}.$$
The function $T^m$
characterizes $m$ up to translation and scaling and the pair $(L,
T^m)$ characterizes $m$ up to translation. In the following discussion, functions will depend on time $t\in [0,1]$ and normalized arc length $s$, also in $[0,1]$. For clarity, we will let $\Om = [0,1]$ for the arc length, i.e., write $t\in [0,1]$, $s\in\Om$. Given a function $h(t,s)$, $(t, s)\in  [0,1]\times \Om$, we will denote $\dot h = \prt_t h$ for derivatives with respect to the time variable, and $dh$ for derivatives with respect to the parameter.\\

We consider the group of diffeomorphisms $G$ of $\Om$, which acts on the set $M$ of measurable functions $a:\Om \to S^{d-1}$ (the unit sphere in $\mR^d$) by $\phi \cdot a =
a\circ \phi^{-1}$ (and the group product is $\phi\psi = \phi\circ \psi$).  ``Tangent vectors'' to $G$ at $\phi = \id$ are functions $v:\Om \to \mR$, with $v(0) = v(1) = 0$, and ``tangent vectors'' to $M$ at $\al$ are functions $\xi: \Om\to \mR^d$ such that $\xi(s)^T\al(s) = 0$ almost everywhere on $\Om$ (we will not attempt here to define $G$ and $M$ as manifolds and rigorously describe their tangent spaces). 

We consider a Hilbert space $V$ of continuous functions $v:\Om\to \mR$ (satisfying $v(0) = v(1) = 0$), with norm denoted $\|\cdot\|_V$ and 
define
\[
\|(v, \xi)\|_{(\id, \al)}^2  = \|v\|_V^2 + \frac{1}{\sig^2} \|\xi\|_2^2
\]
where $\|\cdot\|_2$ is the $L^2$ norm. Then, the metamorphosis energy in \eqref{eq:meta.gen} becomes
\begin{equation}
\label{eq:meta.1}
\int_0^1 \|\dot \phi \circ \phi^{-1}(t)\|_V^2 \, dt + \frac{1}{\sig^2} \int_0^1 \int_\Om |\dot\al \circ \phi^{-1}|^2 \, dt,
\end{equation}
to be minimized with $\al(0) = a_0$ and $\al(1) = a_1\circ \phi(1)$.

When the Hilbert space $V$ is continuously embedded in $C^1(\Om, \mR)$, and the functions $\al$ are differentiable at all times, this formulation (that can be referred to as a Lagrangian form of metamorphosis) have an equivalent Eulerian form obtained by letting $v = \dot \phi \circ \phi^{-1}$ and $a = \al \circ \phi^{-1}$, for which \eqref{eq:meta.1} can be rewritten as
\begin{equation}
\label{eq:meta.2}
\int_0^1 \|v(t)\|_V^2 \, dt + \frac{1}{\sig^2} \int_0^1 \int_\Om |\dot a(t) + v(t) da(t)|^2 \, dt
\end{equation}
minimized with $a(0) = a_0$ and $a(1) = a_1$. The equivalence comes from the fact that the equation $\dot \phi(t) = v(t)\circ\phi(t)$ has a unique solution when $v$ is $C^1$. This is the standard form of metamorphosis discussed in  \cite{ty05b,hty09}. 
\\

In our case, however, functions in $V$ are not necessarily differentiable, and we will work only from the first formulation, \eqref{eq:meta.1}. More precisely, we take
\begin{equation}
\label{eq:n.v.1}
\|v\|_V^2 = \int_\Om dv(s)^2\, ds,
\end{equation}
which implies that functions $v\in V$ are continuous and satisfy a H\"older condition of order $q$ for any $q<1/2$, but are not necessarily Lipschitz continuous. This choice being made, an elementary computation followed by a change of variable in both integrals provides our final expression of the metamorphosis energy, namely
\begin{equation}
\label{eq:U.sig}
U_\sig(\phi, \al) = \int_0^1\int_\Om \frac{( d\dot \phi)^2}{d\phi}\, ds\,dt + \frac1{\sig^2}\int_0^1\int_\Om |\dot \al|^2 d\phi\, ds\,dt
\end{equation}
which needs to be minimized over all trajectories $t\mapsto \phi(t)$ and $t\mapsto \al(t)$, such that $\phi$ is at all times an increasing diffeomorphism of $\Om$ and $\al$ a function from $\Om\to S^{d-1}$, with boundary conditions $\al(0) = a_0$ and $\al(1) = a_1\circ \phi(1)$. This energy coincides (up to a multiplicative constant) with the one introduced in \cite{mio2007shape}.

\subsection{First Reduction}
We consider the minimization of $U_\sig$ with respect to $\al(\cdot)$, with given $\phi(\cdot)$ and use this to reduce the original $d$-dimensional problem to a similar two-dimensional one.  Indeed, we first notice that, in order to minimize the second term in $U$, it suffices to minimize separately each integral
\begin{equation}
\label{eq:meta.data.cost}
\int_0^1 |\dot\al(t,s)|^2 d\phi(t,s) dt
\end{equation}
where $s$ is fixed. 
Considering this integral, we write $\al(t,s) = \tilde \al(\la(t,s), s)$ where
\[
\la(t,s) := \frac{\int_0^t d\phi(t,s)^{-1}\,dt}{\int_0^1 d\phi(t,s)^{-1}\,dt}
\]
is an increasing function satisfying $\la(0,s) = 0$ and $\la(1,s) = 1$.
We have $\dot \al = \dot \la \dot{\tilde\al}(\la, s)$ and 
\[
\int_0^1 |\dot\al(t,s)|^2 d\phi(t,s) dt = \frac{1}{c(s)}\int_0^1 |\dot {\tilde\al} (\la(t,s), s)|^2 \dot \la(t,s)\, dt = \frac{1}{c(s)} \int_0^1 |\dot {\tilde\al} (t, s)|^2\, dt
\]
with 
\[
c(s) = \int_0^1 d\phi(t,s)^{-1}\,dt.
\]
This integral must be minimized subject to $\tilde \al(0,s) = a_0(s)$, $\tilde \al(1,s) = a_1\circ \phi_1(s)$ and $|\tilde \al(t,s)| = 1$ for all $t$, and the solution is given by  the arc of circle between $\tilde\al(0,s)$ and $\tilde\al(1,s)$, which can be expressed as
\[
\tilde \al(t, s) = \frac{\sin((1-t) \om(s))}{\sin\om(s)} a_0(s) + \frac{\sin(t \om(s))}{\sin\om(s)} a_1\circ \phi_1(s)
\]
where
\[
\om(s) = \mathrm{arccos} (a_0(s)^T a_1\circ \phi_1(s)).
\]
The optimal $\al$ is therefore given by 
\begin{equation}
\label{eq:opt.alpha}
\al(t, s) = \frac{\sin((1-\la^{-1}(t,s)) \om(s))}{\sin\om(s)} a_0(s) + \frac{\sin(\la^{-1}(t,s) \om(s))}{\sin\om(s)} a_1\circ \phi_1(s)
\end{equation}
where $\la^{-1}(t,s)$ is defined by $\la(\la^{-1}(t,s),s) = t$. The optimal cost in \eqref{eq:meta.data.cost} is then $\om(s)^2/c(s)$. We notice that, because $\om(s)\in [0, \pi]$, the coefficients in \eqref{eq:opt.alpha} are non-negative. Moreover, $\al(t,s)$ is at all times in the plane generated by $a_0(s)$ and $\al_1 = a_1\circ \phi_1(s)$. \\

We now fix $\phi_1$ and study optimal metamorphoses with $\phi(1, \cdot) = \phi_1$. Introduce the vector $a_0^\perp$ perpendicular to $a_0$ in the plane generated 
by $a_0$ and $\al_1$, defined by
\[
\al_1 = \cos\om \, a_0 + \sin\om\, a_0^\perp.
\]
(This is well defined if $\om \in (0, \pi)$ and we choose $a_0^\perp$ arbitrarily otherwise.)
Without loss of generality, we can search for optimal metamorphoses taking the form
\[
\al(t,s) = \cos\tau(t,s) a_0(s) + \sin\tau(t,s) a_0^\perp(s)\,.
\]
Letting $\xi(t, s) = (\cos\tau(t,s), \sin\tau(t,s))\in S^1$, we can write
$U_\sig(\phi, \al) = \tilde U_\sig(\phi, \xi)$, where 
\begin{equation}
\label{eq:U.sig.th}
\tilde U_\sig(\phi, \xi) = \int_0^1\int_\Om \frac{( d\dot \phi)^2}{d\phi}\, ds\,dt + \frac1{\sig^2}\int_0^1\int_\Om |\dot \xi|^2 d\phi\, ds\,dt.
\end{equation}
This function now has to be minimized  subject to $\phi(0, \cdot)  = \id$, $\phi(1, \cdot) = \phi_1$, $\xi(0, \cdot) = (1,0)$ and $\xi(1, \cdot)  = (\cos \om, \sin\om)$, with $\om = \mathrm{arccos} (a_0^T a_1\circ \phi(1, \cdot))$. In other terms, we have reduced the $S^{d-1}$-valued metamorphosis problem to an $S^1$-valued problem, or, equivalently, our metric on $d$-dimensional curves to a two-dimensional case.
%Note that we get an identical problem if we replace the coundary conditions by  $\th(0, \cdot) = \th_0 := 2\ell(s)\pi$ and $\th(1, \cdot) = \th_1:= \om + 2\ell(s)\pi$, where $\ell(s)$ is anu measurable interger-valued function.

%We now consider the problem in which $\phi(1, \cdot) = \phi_1$ is fixed (optimizing later for the correspondence $\phi_1$) and define a new variable with respect to which $\tilde U$ will take a simpler form. Note that fixing $\phi_1$ also fixes the function $\om$.
%For fixed $s$, the function $t \mapsto \al(t, s)$ is continuous and can be lifted to take the form $\al(t, s) =  (\cos \th(t, s), \sin\th(t, s))$ where $\th(\cdot, s)$ is a continuous function defined on $[0, 1]$. The function $\th$ is uniquely defined as soon as one chooses a representation $a_0(s) = (\cos \th_0(s), \sin\th_0(s))$, which we assume to be given from now on. Notice that, with this representation, we have
%\[
%|\dot \al(t, s)|^2 = \frac{\dot \th(t, s)^2}{4\sig^2} 
%\] 
%($\th$ obviously inherits the time differentiability of $\al$).

\subsection{Second Reduction}
The second reduction is the by-now well known square root transform that will move the problem into a standard Hilbert framework.
 Because $\xi(t,s)$ is differentiable in time, one can define uniquely a differentiable function $\tau(t,s)$ such that $\xi(0,s) = (1,0)$ and $\xi(t,s) = (\cos\tau(t,s), \sin\tau(t,s))$ at all times. 
Define $q(t,s)$ by
\begin{equation}
\label{eq:sqvf}
q(t, s) =\sqrt{d \phi(t, s)} \left(\cos \eta(t,s), \sin\eta(t,s)\right).
\end{equation}
with $2\sig\eta(t, s) = \tau(t, s)$.
%where $\ep: \Om \to \{-1, 1\}$ is measurable. 
Then, a straightforward computation yields 
\[
4|\dot q(t, s)|^2 = \frac{d\dot\phi^2}{d \phi} + \frac{1}{\sig^2} |\dot\xi(t, s)|^2
\]
so that
\begin{equation}
\label{eq:U.q}
\tilde U_\sig(\phi, \xi) = 4\int_0^1 \|\dot q(t, \cdot)\|^2_2\, dt
\end{equation}
We also note that 
\[
\|q(t, \cdot)\|_2^2  = \int_\Om d\phi(t, s)\, ds = 1
\]
so that $q(t, \cdot)$ is a curve on the unit sphere of $L^2(\Om, \mR^2)$. This implies that its energy,
$\int_0^1 \|\dot q(t, \cdot)\|^2_2\, dt$, cannot be larger than that of the minimizing geodesic on this unit sphere, which is the shortest great circle connecting the  functions $q_0 = \left(1, 0\right)$ (which is constant) 
and $q_1 = \sqrt{d\phi_1} \left(\cos\left(\eta(1,\cdot)\right), \sin\left(\eta(1, \cdot \right)\right)$.
% with the usual notation $\th(1) = \th(1, \cdot)$.
Letting
 \[
 \rho = \mathrm{arccos} \scp{q_0}{q_1}_2 = \mathrm{arcos} \int_\Om \sqrt{d\phi_1}\cos\left(\frac{\om(s)}{2\sig}\right)\, ds, 
 \]
 this geodesic is given by 
\begin{equation}
\label{eq:great.circ} 
\ga(t, s) = \frac{\sin((1-t)\rho)}{\sin\rho} q_0(s) + \frac{\sin(t\rho)}{\sin\rho} q_1(s)
\end{equation}
with energy equal to $\rho^2$. We therefore find that 
\begin{equation}
\label{eq:lb.U}
\tilde U_\sig(\phi, \xi) \geq 4\, \mathrm{arcos}^2 \int_\Om \sqrt{d\phi_1(s)}\cos\left(\frac{\om(s)}{2\sig}\right)\, ds\,.
\end{equation}
\\

This provides a lower-bound for the metamorphosis energy. To prove that this lower-bound is achieved, we now investigate whether  the trajectory $\ga$ in \eqref{eq:great.circ} can be derived from a valid trajectory $(\psi(\cdot, \cdot), \mu(\cdot, \cdot))\in G\times M$ that connects $(\id, 0)$ to $(\phi_1, \xi_1)$. 
%Assume that a trajectory $q$ on the unit sphere of $L^2(\mR^2)$ is given. We look for representations

We are therefore looking for representations of $\ga$ in the form
\[
\ga(t, s) =\sqrt{d \psi(t, s)} \left(\cos\tilde \eta(t, s), \sin\tilde\eta(t, s)\right),
\]
with  $\tilde \eta(0, s) = 0$, which uniquely defines $\tilde \eta(t, s)$ by continuity in $t$. Notice that we automatically have $d\psi(0, \cdot) = 1$ and $d\psi(1, \cdot) = d\phi_1$ by definition of $q_0$ and $q_1$. For $t\in(0,1)$, we have
\[
d \psi(t, s) = |\ga(t, s)|^2,
\]
so that $\psi(t, \cdot)$ is non-decreasing and satisfies $\psi(t,0) = 0$, $\psi(t, 1) = 1$. The function $s\mapsto \dot \psi(t,s)$  is positive if and only if $q(t, s)$ does not vanish, which requires 
\begin{multline}
\label{eq:not.pi0}
\sin^2((1-t)\rho)  + \sin^2(t\rho) d\phi(1) \\
+ 2\sin((1-t)\rho) \sin(t\rho) \sqrt{d\phi_1}\cos\left(\frac{\om}{2\sig}\right) >0.
\end{multline}
A sufficient condition for this to holds for all $t$ is that the cosine term is strictly larger than $-1$, 
which is equivalent to 
\begin{equation}
\label{eq:not.2pi}
\frac{\om}{2\sig} \not \equiv \pi \quad (2\pi).
\end{equation}
Since $\om \in [0, \pi]$, this condition will be automatically satisfied if $2\sig > 1$.\\
%We will return to this condition (which means that the left-hand term is not a multiple of $k\pi$ with odd $k$) later in this discussion.
%It it is not true, then the left-hand side of \eqref{eq:not.pi0} is $(\sin((1-t)\rho) - \sqrt{d\phi_1} \sin(t\rho))^2$, and this term vanishes for some $t\in (0,1)$.   

Because $\ga(1) = q_1$, we must have
\begin{equation}
\label{eq:2sig.int}
\tilde\eta(1, s) = \om(s)/2\sig + 2 k(s) \pi
\end{equation}
where $k$ is an integer-valued function. The right-hand side of  \eqref{eq:great.circ} is a linear combination of $q_0$ and $q_1$ with positive coefficients, which implies that the time-continuous angular representation of $q(t)$ starting at $0$ cannot deviate by more than $\pi$ from its initial value, i.e.,
\begin{equation}
\label{eq:2sig.1}
-\pi \leq \tilde\eta(1,s) \leq \pi,
\end{equation}
%For the same reason, using equation \eqref{eq:opt.alpha}, we have
and we also have
\begin{equation}
\label{eq:2sig.2}
0 \leq \om(s) \leq \pi.
\end{equation}
We now (and for the rest of the discussion) make the assumption that  $2\sig \geq 1$. Under this assumption,  $\tilde\eta(1, \cdot) = \om/2\sig$ satisfies  both \eqref{eq:2sig.int} and \eqref{eq:2sig.1}. Since it is clear that only one value of $\tilde\eta(1,s)$ can satisfy the two equations together, we find that, for all $s\in \Om$, one has $\tilde\eta(1,s) = \om(s)/\sig$, and the curve $\ga$ is associated with a trajectory $(\psi, \be)$ between $(\id, 0)$ and $(\phi_1, \xi_1)$.  \\

We have therefore proved that $\ga$ in \eqref{eq:great.circ} provides a valid improved solution to  the original $(\phi, \al)$ as soon as \eqref{eq:not.2pi} is satisfied, which is true as soon as $2\sig > 1$. 

If $\sig =1/2$, then $\eqref{eq:not.2pi}$ may not hold for the curve in \eqref{eq:great.circ}. However, the minimum of $U_\sig(\phi, \al)$ with given $\phi(1) = \phi_1$ is still given by the geodesic energy of this curve. To see this, it suffices to  consider a small variation $\tilde a_0$ of $a_0$ such that $\tilde a_0^T a_1\circ \phi_1 > -1$, so that  \eqref{eq:not.2pi} is satisfied  with $\tilde a_0$ instead of $a_0$, and the minimum energy when starting from $\tilde a_0$ is the geodesic energy of the associated great circle. One can then use the fact that $U$ is a geodesic energy for a Riemannian metric on $G\times M$, and combine this with the triangular inequality  for the sequence of  geodesics going from  $(\id, a_0)$ to $(\id, \tilde a_0)$ then to $(\phi(1), a_1)$. Indeed, the energy of the sequence is larger than the minimal energy between $(\id, a_0)$ and $(\phi(1), a_1)$, but arbitrarily close to the energy of the minimal geodesic between $(\id, \tilde a_0)$ then to $(\phi(1), a_1)$, itself arbitrarily close to the lower bound in \eqref{eq:lb.U}.\\

We summarize this discussion in  the following theorem.
\begin{theorem}
\label{th:meta.curve.sig}
Assume that $2\sig\geq 1$, and let $\phi_1:\Om \to\Om$ satisfy $\phi_1(0) = 0$, $\phi_1(1) = 1$ and $d\phi_1 >0$. Then
\begin{multline}
\label{eq:meta.curve.sig}
\inf\defset{U_\sig(\phi, \al): \phi(1) = \phi_1, \al(0) = a_0, \al(1) = a_1\circ \phi_1} \\
= 4\, \mathrm{arccos}^2 \int_\Om \sqrt{d\phi_1(s)}  \cos\left(\frac{\mathrm{arccos} (a_0(s)^Ta_1\circ\phi_1(s))}{2\sig}\right) \, ds.
\end{multline}
Moreover, if $2\sig >  1$, the minimum is achieved and can be deduced from for a geodesic curve $\ga$ on the unit sphere of $L^2(\mR)$. 

This induces a distance on $M$, given by
\begin{multline}
\label{eq:meta.curve.sig.dist}
d_\sig(a_0, a_1)  = 2 \inf_{\phi_1}\left(\mathrm{arccos} \int_\Om  \sqrt{d\phi_1(s)}\cos\left({\mathrm{arccos} (a_0(s)^Ta_1\circ\phi_1(s))}/{2\sig}\right)\,ds\right)
\end{multline}
minimized over all strictly increasing diffeomorphisms of $\Om$.

\end{theorem}
This theorem is essentially proved in the discussion that precedes it, in which we have left a few loose ends, mostly regarding measurability and dealing with sets of measure 0 that can be  tied without too much effort by an interested reader. 

It is also interesting to express the distance in terms of the original curves, say $m_0$ and $m_1$, of which $a_0$ and $a_1$ are the unit tangents. Indeed, if $m$ is a parametrized curve (not necessarily with arc length), then (still assuming unit length), its associated tangent function is $a_m = T^m \circ s_m^{-1}$ where $T_m = dm/|dm|$ is the unit tangent in the original parametrization and $s_m$, such that $ds_m = |dm|$, is the arc length reparametrization. The distance is therefore given by
\[
d_\sig(m_0, m_1)  = 2 \inf_{\phi_1}\left(\mathrm{arccos} \int_\Om  \sqrt{d\phi_1} F_\sig(a_{m_1} \circ \phi_1, a_{m_0})\,ds\right)
\]
where we have denoted, for short
\[F_\sig(a_1, a_0) = \cos\left({\mathrm{arccos} (a_0(s)^Ta_1)}/{2\sig}\right).
\]
We therefore have
\begin{align*}
d_\sig(m_0, m_1)  &= 2 \inf_{\phi_1}\left(\mathrm{arccos} \int_\Om  \sqrt{d\phi_1} F_\sig(T^{m_1}\circ s_{m_1} \circ \phi_1, T^{m_0}\circ s_{m_0})\,ds\right)\\
&= 2 \inf_{\phi_1}\left(\mathrm{arccos} \int_\Om  |dm_0| \sqrt{d\phi_1\circ s_{m_0}} F_\sig(T^{m_1}\circ s^{-1}_{m_1} \circ \phi_1\circ s_{m_0}, T^{m_0})\,du\right)
\end{align*}
If we let $\psi_1 = s_{m_1}^{-1}\circ \phi_1 \circ s_{m_0}$, we get the alternative expression
\begin{multline}
\label{eq:dist.m}
d_\sig(m_0, m_1)  = \\
2 \inf_{\psi_1}\left(\mathrm{arccos} \int_\Om  \sqrt{\psi_1} \sqrt{|dm_1|}\circ \psi_1 \sqrt{|dm_0|}F_\sig(T^{m_1}\circ \psi_1, T^{m_0})\,du\right),
\end{multline}
which expresses the distance directly in terms of the compared curves.

\subsection{The Square Root Velocity Function.} In the two-dimensional case, one can represent functions $\al: \Om\to S^1$ in the form $(\cos\th, \sin\th)$ for some angle function $\th$, defined up to the addition of a multiple of $2\pi$. Given such a representation, one can then consider the transform 
\[
\mathcal G: (\phi, \th) \mapsto q = \sqrt{d\phi}\, (\cos\th/2\sig, \sin\th/2\sig)
\]
that defines a mapping from $\Diff(\Om) \times L^2(\Om, \mR)$ to the unit sphere of $L^2(\Om, \mR^2)$. Adding a time dependency, we find, using this transform, that 
\[
\int_0^1\int_\Om \frac{( d\dot \phi)^2}{d\phi}\, ds\,dt + \frac1{\sig^2}\int_0^1\int_\Om |\dot \xi|^2 d\phi\, ds\,dt = 4\int_0^1 \|\dot q\|_2\, dt, 
\]
so that minimizers on the left can be associated with geodesics on the unit sphere. One can then compute the metamorphosis distance by minimizing the lengths of great circles between, say,  $\mathcal G(\id, \th_0)$ and $\mathcal G(\phi_1, \th_1\circ \phi_1)$, for a given $\phi_1$, and optimizing over all angle representations $\th_0, \th_1$ of $a_0$, $a_1$ that satisfy the constraint
\[
-2\sig\pi \leq \th_1 \circ \phi_1 - \th_0 \leq 2\sig\pi\,,
\]
because \eqref{eq:2sig.1} still needs to hold for any time-continuous angle representation of $q$. 

This provides the same distance $d_\sig$ as the one obtained in Theorem \ref{th:meta.curve.sig}. Notice, however, that this construction is special to the two-dimensional case. In dimension $d>2$, our reduction to unit-sphere geodesics depended on the end-points $a_0$ and $a_1$, and could not be deduced from a direct transformation applied to curves themselves, such a $\mathcal G$. The only exception is the case $\sig=1/2$, for which $\mathcal G$ is equivalent to $(\phi,\al) \mapsto \sqrt{d\phi} \,\al$. This transform, called the ``square root velocity function,'' is clearly applicable to arbitrary dimensions. It  has been extensively studied in the literature, and we refer to \cite{srivastava2016functional} and references within for additional details and applications.\\

Returning to the two-dimensional case, alternate expressions of the distance can be derived for simple values of $\sig$ given angle representations $\th_0$ and $\th_1$ for $a_0$ and $a_1$. Indeed,  the cosine in \eqref{eq:meta.curve.sig.dist} is  given by $\cos(\th_1\circ \phi_1 - \th_0)$ if $\sig=1/2$, by $|\cos((\th_1\circ \phi_1 - \th_0)/2)|$ if $\sig=1$, and by
\[
\max(|\cos((\th_1\circ \phi_1 - \th_0)/4)|, |\sin((\th_1\circ \phi_1 - \th_0)/4)|)
\]
if $\sig=2$.\\

The case $\sig = 1$ for plane curves was first investigated in \cite{you96,you98}. It has interesting additional features,  because it corresponds to a Riemannian metric on spaces of curves associated with the first-order Sobolev norm of vector fields along the curves  \cite{ymsm07} (see section \ref{sec:closed}). The optimization of $\phi_1$ (which is still needed to compute the distance) can be done efficiently by dynamic programming, and the reader is referred to \cite{you99,ty00} for more details.  \\

To conclude this section, we notice that the metamorphosis metric is invariant to the action of rotations, so that one can optimize  a rotation parameter in all cases considered above. The rotation invariant version of the distance between plane curves for $\sig=1$, for example, is
\begin{multline}
\label{eq:meta.curve.sig.1.dist}
d_{1, \mathrm{rot}} (a_0, a_1) = 2 \inf_{\phi_1,c} \left( \mathrm{arccos} \int_\Om \sqrt{d\phi_1} \left|\cos\left(\frac{\th_1\circ \phi_1(s) - \th_0(s) -c}{2}\right)\right| \, ds\right)
\end{multline}
where $c$ is a scalar. In higher dimensions, one needs to optimize \eqref{eq:meta.curve.sig.dist} with $a_0$ replaced by $Ra_0$ when $R$ varies over all rotations of the $d$-dimensional space. 
%Using formula \eqref{eq:meta.curve.sig.dist}, the distance is  minimizes, for $\sig = 1$

%The case $d=1$ is also interesting, in which a function $s\mapsto m(s)$  

\subsection{One-dimensional case}

Curves in one dimension are functions $u \mapsto f(u)$ and the unit tangent is $T^f(u) = \sign(df(u))$, assuming that the latter is non-zero almost everywhere. This reduces the previous representation  to functions $a:\Om \to \{-1, 1\}$, which does not leave much room for the definition of time-continuous metamorphoses, so that this approach cannot be directly extended to this case. 

One can however bypass this difficulty by  associating plane curves to such functions. Defining $\sign(0) = 0$, we can associate to a function $f: [0,1]\to\mR$, the horizontal curve $m_f(u)$ such that $m_f(0) = 0$ and 
\[
dm_f(u) = (|df(u)| \sign(df(u)), 0).
\]
  The normalization to length one means that we normalize function by their total variation
\[
TV(f) = \int_\Om |df(u)|\,du
\]
and the arc length is $s = \psi_f(u)$ with $d\psi_f = |df|$. To simplify the discussion, we will assume that $df$ vanishes over no interval so that $\psi_f$ is strictly increasing. The ``angle function'' associated with $m_f$ is therefore $a_f(s) = \sign(df)\circ \psi_f^{-1}$. Given two functions $f_0$ and $f_1$, we use the plane curve distance to compare $a_{f_0}$ and $a_{f_1}$, noting that in this case, one has $\om(s) =\pi$ if $a_{f_1}\circ\phi_1(s) \neq a_{f_0}(s)$ and $0$ otherwise. One therefore gets
\begin{multline}
\label{eq:meta.curve.sig.dist.1D}
d_\sig(f_0, f_1)  = 2 \inf_{\phi_1}\bigg(\mathrm{arccos} \int_\Om  \sqrt{d\phi_1(s)}\Big(\bfone_{a_{f_0}(s) = a_{f_1}\circ\phi_1(s)} \\
+ \cos\left(\frac{\pi}{2\sig}\right) \bfone_{a_{f_0}(s)\neq a_{f_1}\circ\phi_1(s)}\Big)\,ds\bigg)\,.
\end{multline}
Notice that this compares functions modulo reparametrization, i.e., $a_0$ and $a_0\circ \phi$ are considered as identical for any increasing diffeomorphism of $\Om$. 
%If this creates too much invariance, one can use the representation $\CG$ for the curves $(a,0)$ without optimizing for reparametrization, which provides the distance 
%\begin{multline}
%\label{eq:meta.curve.sig.dist.1D.2}
%\delta_\sig(a_0, a_1)  = 2\mathrm{arccos} \int_\Om  \sqrt{da_0(s)da_1(s)}\Big(\bfone_{a_0(s) = a_1(s)} \\
%+ \cos\left(\frac{\pi}{2\sig}\right) \bfone_{a_0(s)\neq a_1(s)}\Big)\,ds\,.
%\end{multline}

\subsection{The smooth case}
\label{sec:meta.curves.2}

In our discussion so far, we have placed little regularity conditions on the functions $a\in M$ beyond their measurability. The resulting class of curves includes, in particular, polygonal curves, for which $a$ is piecewise constant. We here briefly discuss the changes that need to be made in the discussion when the considered curves are  smooth.

In this case, the integrals in \eqref{eq:meta.data.cost}  cannot be minimized independently for each $s$, because we need to ensure that the solution that one obtains is a continuous value of $s$. However, even if not necessarily a minimizing geodesic on $S^{d-1}$, the function  $t\mapsto \tilde\al(t,s)$ for fixed $s$ must still be locally minimizing, i.e., it must still be supported by the great circle connecting $a_0(s)$ and $a_1\circ \phi_1(s)$. This implies that the optimal solution is  still given by \eqref{eq:opt.alpha} with, this time, $\om$ being a continuous lift of $s\mapsto \mathrm{arccos} (a_0(s)^T a_1\circ \phi_1(s))$. One can therefore still reduce the $d$-dimensional setting to a two-dimensional one. 

Taking the same definition for $q$, we find that the metamorphosis energy is no larger than four times the geodesic energy of $q$ in the unit sphere of $L^2(\Om, \mR^2)$. When proving that the lower-bound is achieved, one finds that the function $k(s)$ in \eqref{eq:2sig.int} must be continuous, hence constant. Here, we can use the fact that one can take $\om(0)\in [0, \pi]$ and use the same argument as in the non-smooth case for $s=0$, yielding $k(0) = 0$ and therefore $k(s) = 0$ for all $s$ since $k$ is constant. However, and regardless of the value of $\sig$, one cannot ensure that \eqref{eq:not.pi0} is satisfied unless the compared curves are close enough (so that their angles are at distance less than $2\sig\pi$ after registration). When computed between curves that are too far apart, curves in $M$ deduced from geodesics on the sphere  will typically develop singularities (and therefore step out of $M$ if this space is restricted to smooth curves). 

%For example, \eqref{eq:opt.alpha}, may switch between geodesics on the unit circle when $\om(s)$ reaches $\pi$, and therefore induce jumps in $\al(t,s)$. If one restricts to smooth trajectories, then $\om(s)$ must be replaced by $\th_1(\phi(1,s)) - \th_0(s)$, where $\th_1$ and $\th_0$ are now required to be continuous angle representations of $a_1$ and $a_0$. There is, moreover, no ambiguity on the choice of the integer-valued function $k$ in \eqref{eq:2sig.int}, which needs to be continuous, hence constant, and therefore vanish because $k(0) = 0$. Finally, the optimization over $k^*$ in \eqref{eq:meta.curve.sig} must also be made over constant functions ($k^*$ does not depend on $s$), but theorem \ref{th:meta.curve.sig} holds with this modification. 
%
%For example, in the case $\sig = 1$,  one can only change the sign of the integral globally, and the resulting minimal energy is 
%\begin{multline}
%\label{eq:meta.curve.sig.1.smooth}
%\min\defset{U_1(\phi, \al): \phi(1) = \phi_1, \al(0) = a_0, \al(1) = a_1\circ \phi_1} \\
%= 4\, \mathrm{arccos}^2 \left|\int_\Om  \sqrt{d\phi_1} \cos\left(\frac{\th_1\circ \phi_1(s) - \th_0(s)}{2}\right) \, ds\right|.
%\end{multline}

The smooth case  has also been studied in \cite{bauer2014constructing}, in which a different transform is proposed, leading to a representation of plane curves in a three-dimensional space. More recently, \cite{kurtek2018simplifying} made a study of the smooth case for planar curves with an approach similar to the one we develop here.

%The case in which $2\sig<1$ simplifies in the smooth case, because the function $k$ in \eqref{eq:2sig.int} must also be constant, which shows that the angle function deduced from the end-point of the sphere's geodesic differs from the geodesic one by a constant, i.e., the resulting two plane curves differ by a rotation. If one optimizes the distance with respect to rotation, i.e., over all $\th_1 + c$ where $c\in [0, 2\pi]$, this difference has no impact and the problem can be addressed for arbitrary $\sig$. 

\subsection{Existence of Optimal Metamorphoses}
To complete the computation of the optimal metamorphosis, one must still optimize \eqref{eq:meta.curve.sig} with respect to the final diffeomorphism $\phi_1$. The resulting variational problem is a special case of those studied in \cite{ty01}, which considered the maximization of functionals taking the form
\[
F(\phi) = \int_\Om \sqrt{d\phi} f(s, \phi(s))\, ds
\]
over the set $\text{Hom}^+$ of all strictly increasing functions $\phi:\Om \to \Om$ satisfying $\phi(0) = 0$ and $\phi(1) = 1$, where $f$ is a function defined on $\Om^2$. We let
\[
\De_f = \sup_{\psi\in \mathrm{Hom}^+} \int_\Om \sqrt{d\psi} f(\psi(s), s)\, ds
\]
and define the diagonal band
\[
\Om_c = \left\{(s, s')\in[0,1]^2\ |\ |s-s'| \leq c \right\}\,.
\]
We give without proof the following result, which is a consequence of Theorem 3.1 in
\cite{ty01}.
\begin{theorem}
\label{cont.case}
Assume that $f\geq 0$ is continuous on $\Om^2$ except on a set $G$ that can be decomposed as a union of a finite number of horizontal or vertical segments. Assume also that, for some
\[
c > \sqrt{1- \left(\frac{\De_f}{\|f\|_\infty}\right)^2}\,,
\]
there does not exist any non empty open vertical or horizontal segment
$(a,b)$ such that $(a,b)\sub\Om_c$ and $f_l$ vanishes on $(a,b)$,  where
\[
f_l(x) = \lim_{\de\to 0} \inf_{|y-x|<\de, y\not\in G} f(y)
\]
is the lower semi-continuous relaxation of $f$.

Then there exists $\phi^*\in\text{Hom}^+$ such that $F(\phi^*) = 
\max\{F(\phi),\phi\in\text{Hom}^+\}$. Moreover, if $\phi$ is a
maximizer of $F$, one has, for all
$s\in \Om$, $(\phi(s), s)\in\Om_f$.
\end{theorem}

Intuitively, $f$ vanishing over vertical or horizontal segments allows for either very small or very large values of $d\phi$ at very little cost, resulting in optimal solutions that may have vanishing derivatives or jump discontinuities. In \eqref{eq:meta.curve.sig.dist} (with $\sig = 1/2$), this happens when the tangents of the compared curves are perpendicular.  When $\sig=1$, this happens when $\th_1$ and $\th_0$ are oriented in opposite directions, i.e., their difference is equal to an odd multiple of $\pi$. For $\sig>1$, however, the cosine in \eqref{eq:meta.curve.sig.dist} never vanishes. We also point out that there is no loss of generality in assuming that $f\geq 0$ in the theorem because, if  $f < 0$ on some rectangle, it is easy to check that any trajectory $(s, \phi(s))$ that enters this rectangle can be improved if it is replaced by a trajectory that moves almost horizontally and/or almost vertically within the rectangle, with the new cost converging to 0 over this region. This can be used to show that there is no change in the minimizer if one replaces $f$ by 0 within the rectangle. \\

One can efficiently maximize $F$ by approximating $f$ by a piecewise constant function taking the form
\[
f(s, \tilde s) = \sum_{k=1}^n f_k \bfone_{R_k}
\]
where $R_1, \ldots, R_n$ is a family of rectangles that partition the unit square and 
with $f_{k} \geq 0$, $k=1, \ldots, n$. One can then show that the minimization can be performed over piecewise linear functions $\phi$, which are furthermore linear whenever they cross the interior of a rectangle. The search for the optimal $\phi$ can then be organized as a dynamic program, and run very efficiently (see \cite{ty00,ty01} for details). This method is used in the experiments presented in Figure \ref{fig:meta.curves.2} in which the optimal correspondence is drawn over an image representing the function $\max(f_\sig,0)$, where
\begin{equation}
\label{eq:f.sig}
f_\sig(s, \tilde s) = \cos\left(\frac{\acos\cos(\th_0(s) - \th_1(\tilde s))}{2\sig}\right).
\end{equation}

\begin{figure}
\centering
\includegraphics[trim=2cm 0cm 2cm 0cm,clip,width=0.3\textwidth]{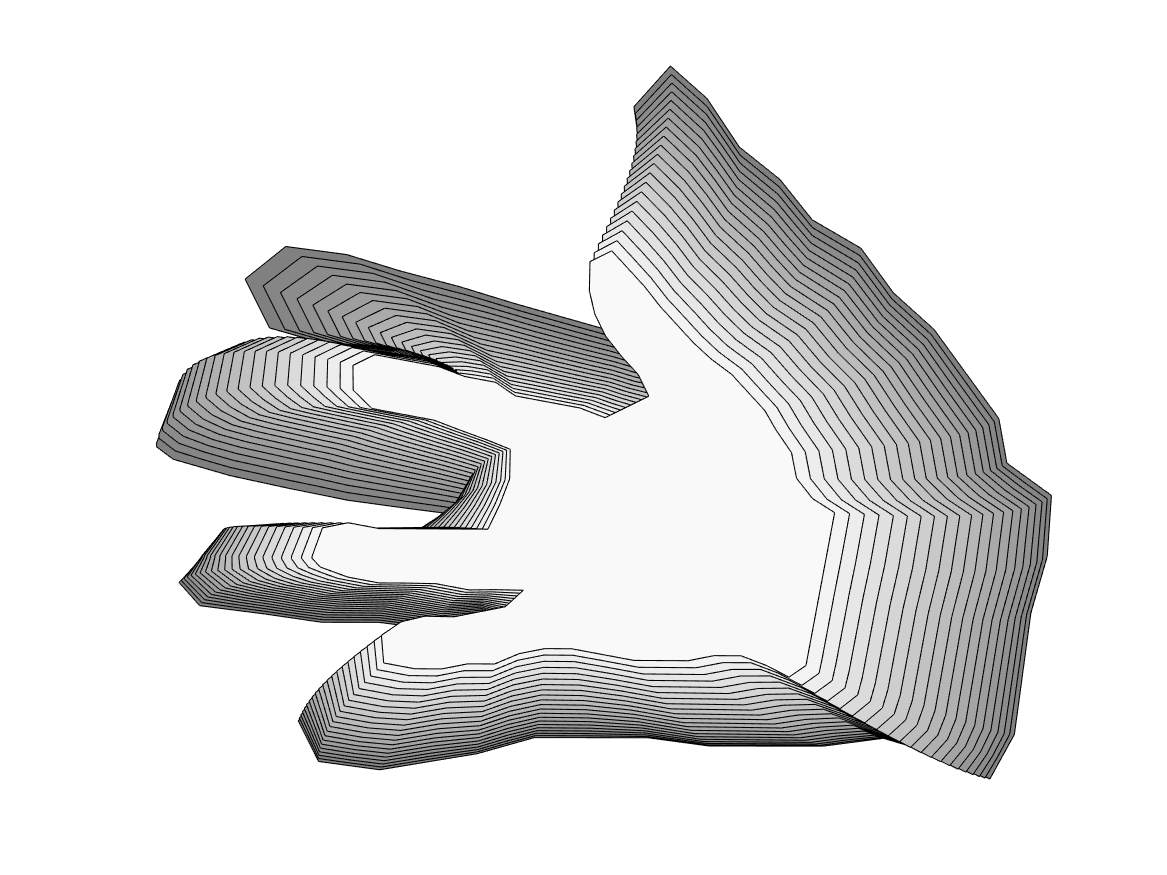}
\includegraphics[trim=2cm 0cm 2cm 0cm,clip,width=0.3\textwidth]{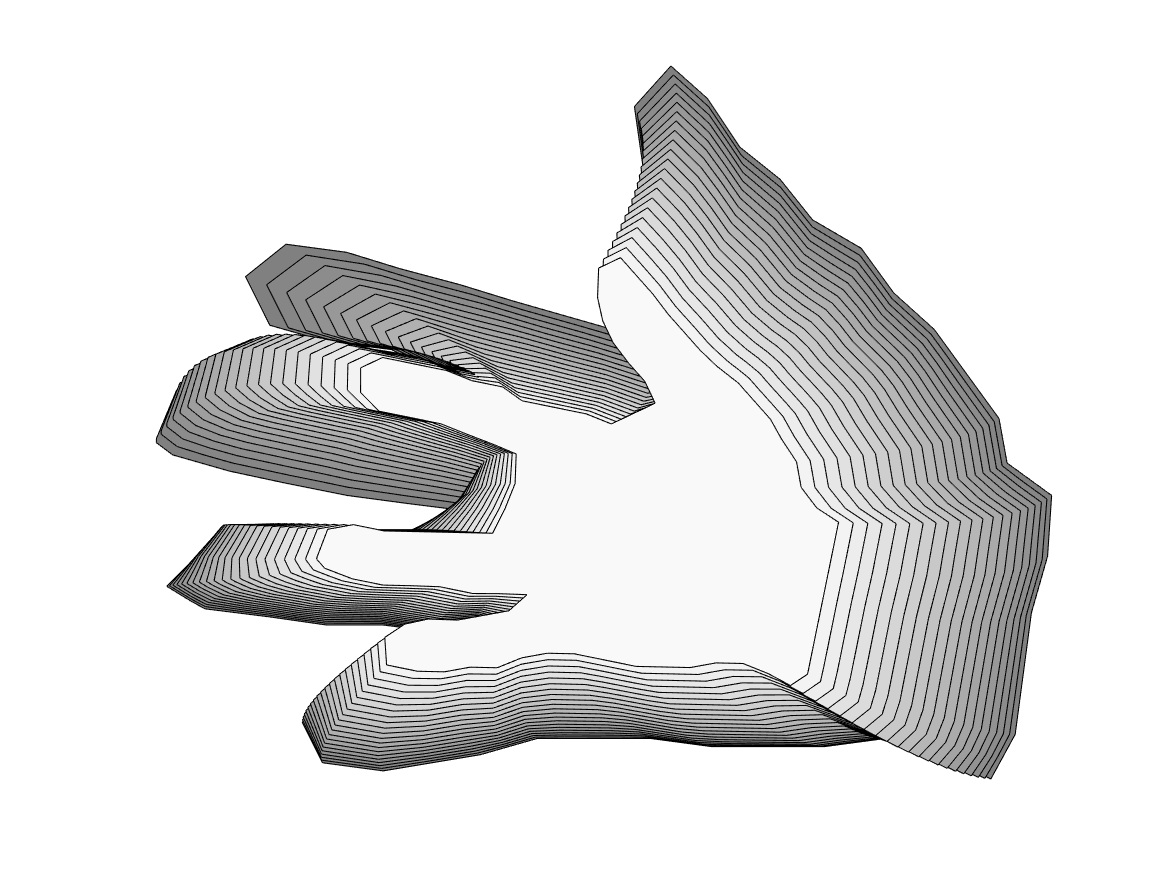}
\includegraphics[trim=2cm 0cm 2cm 0cm,clip,width=0.3\textwidth]{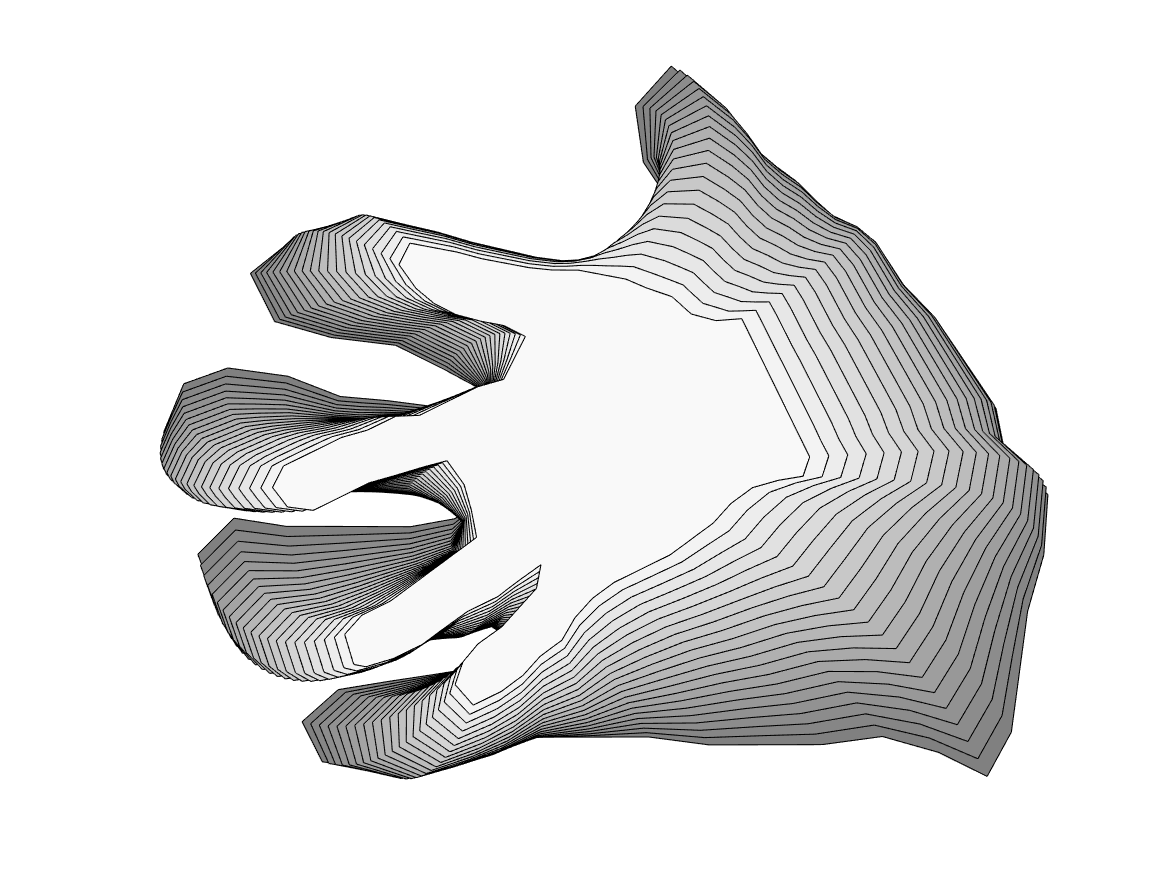}\\
\includegraphics[trim=4cm 0cm 4cm 0cm,clip,width=0.3\textwidth]{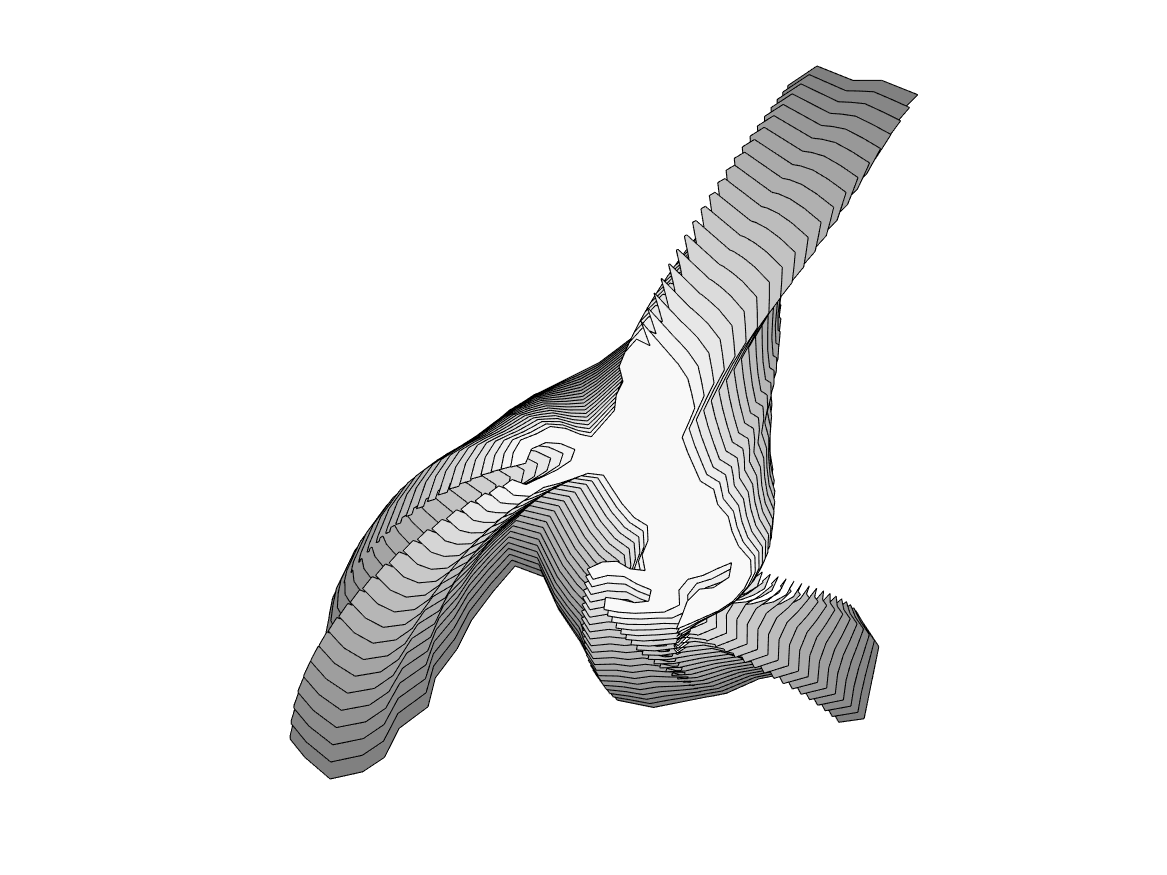}
\includegraphics[trim=4cm 0cm 4cm 0cm,clip,width=0.3\textwidth]{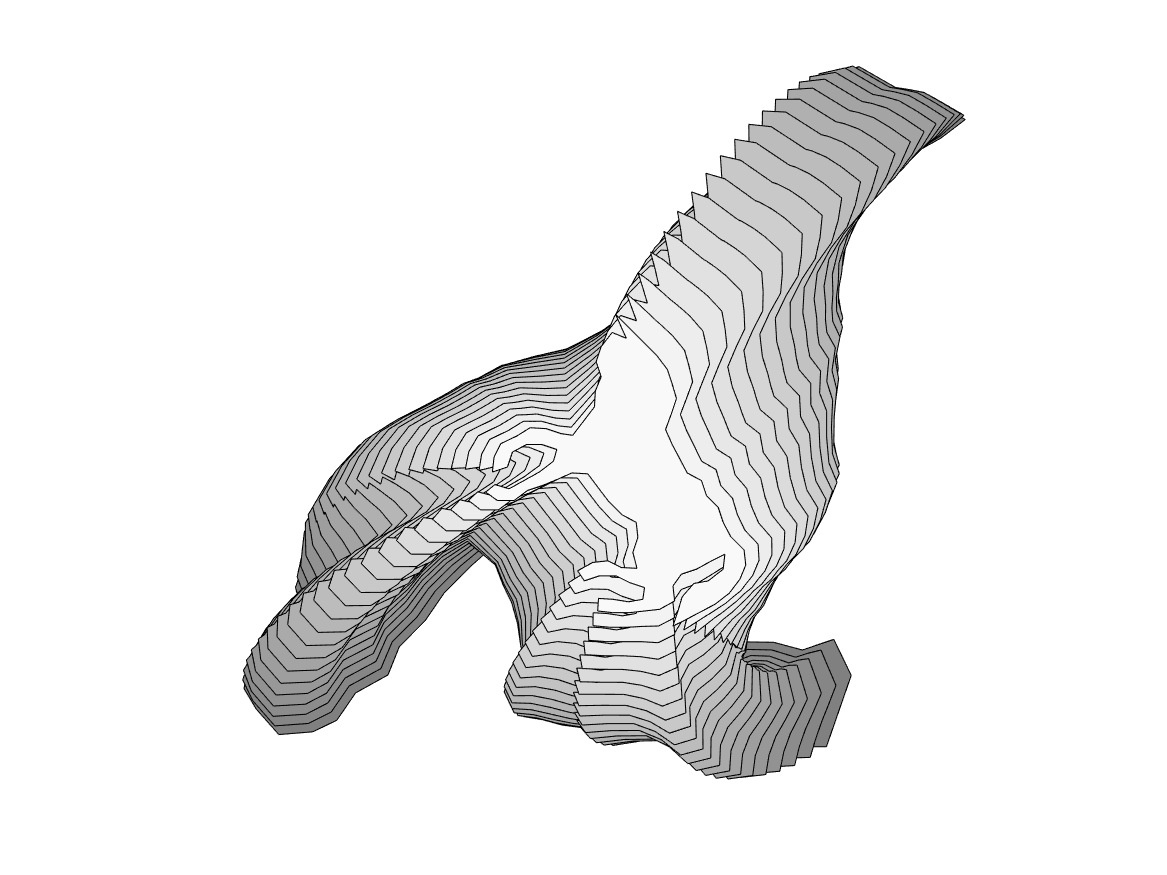}
\includegraphics[trim=4cm 0cm 4cm 0cm,clip,width=0.3\textwidth]{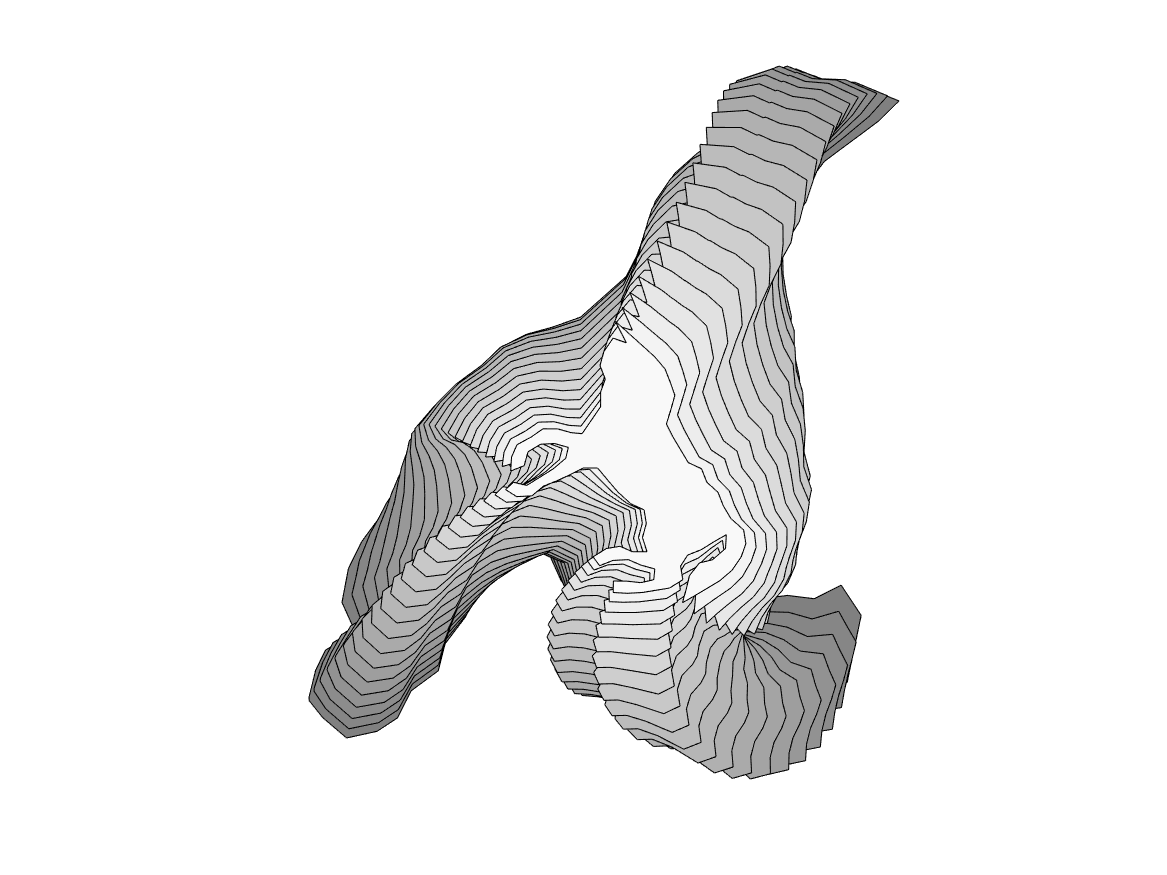}\\
\caption{Two sets of examples of optimal metamorphoses between curves with $\sig = 0.5$ (left), 1 (center) and 2 (right).  Note that, in the first row, the geodesics with $\sig=0.5$ and $\sig=1$ grow the missing finger out of the second finger, which that with $\sig=2$ grows a new thumb.  
In
each image, the geodesic is computed between the outer curve and the
inner curve. Scaling is for visualization only (the computation is
scale invariant)\label{fig:meta.curves}}
\end{figure}

\begin{figure}
\centering
\includegraphics[trim=2cm 0cm 2cm 0cm,clip,width=0.3\textwidth]{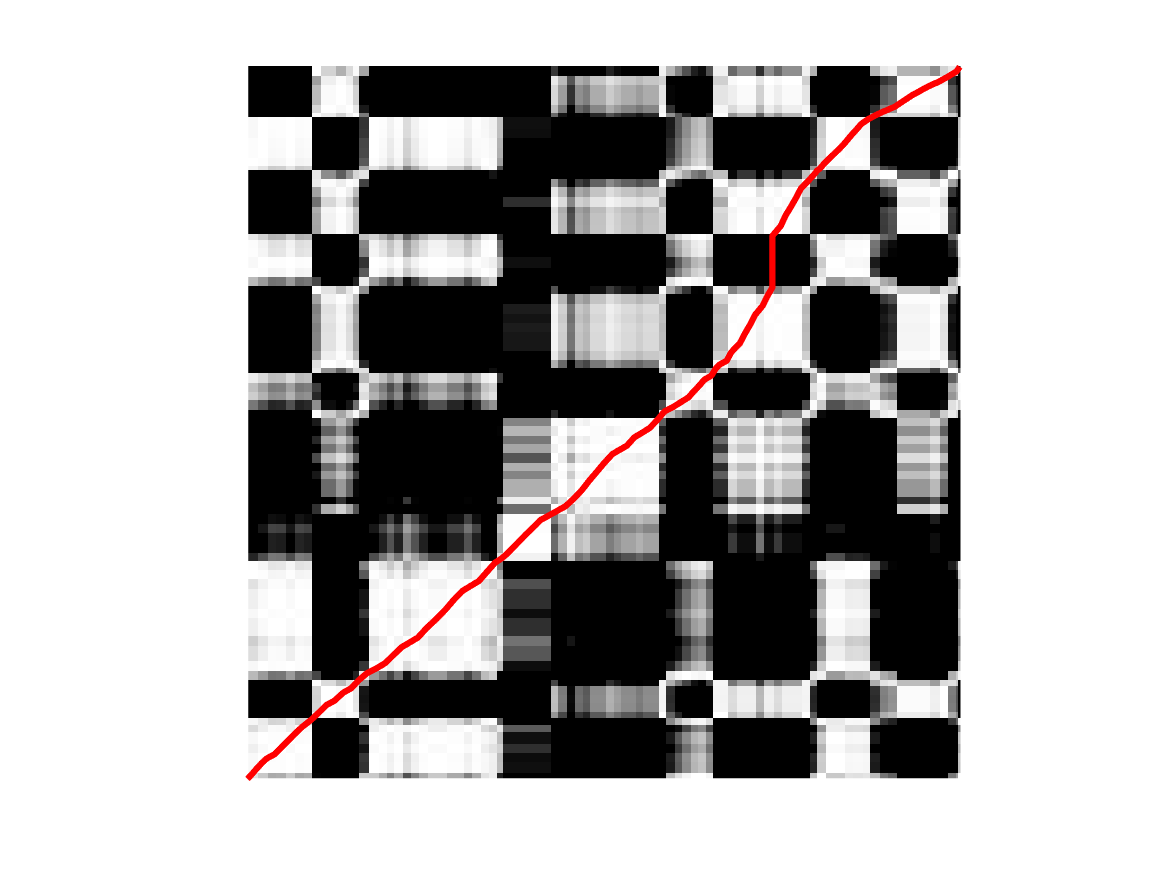}
\includegraphics[trim=2cm 0cm 2cm 0cm,clip,width=0.3\textwidth]{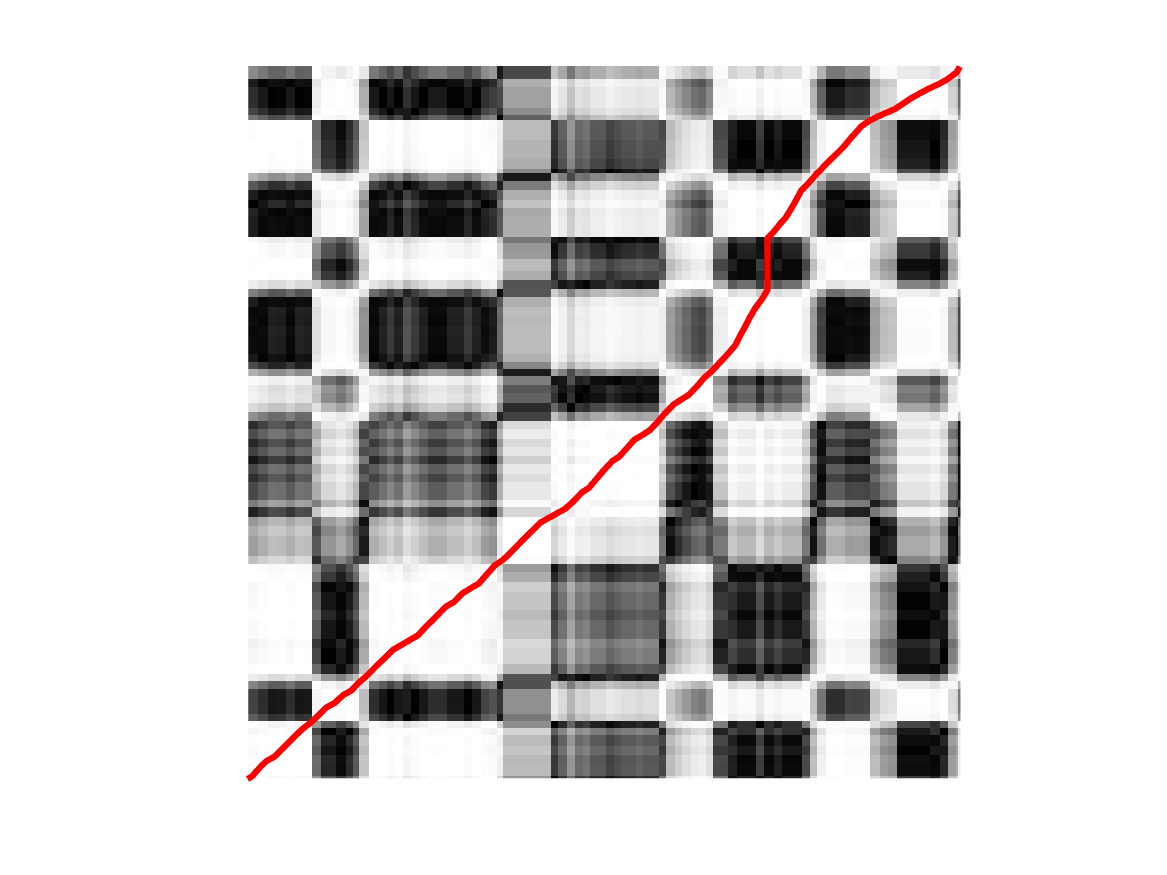}
\includegraphics[trim=2cm 0cm 2cm 0cm,clip,width=0.3\textwidth]{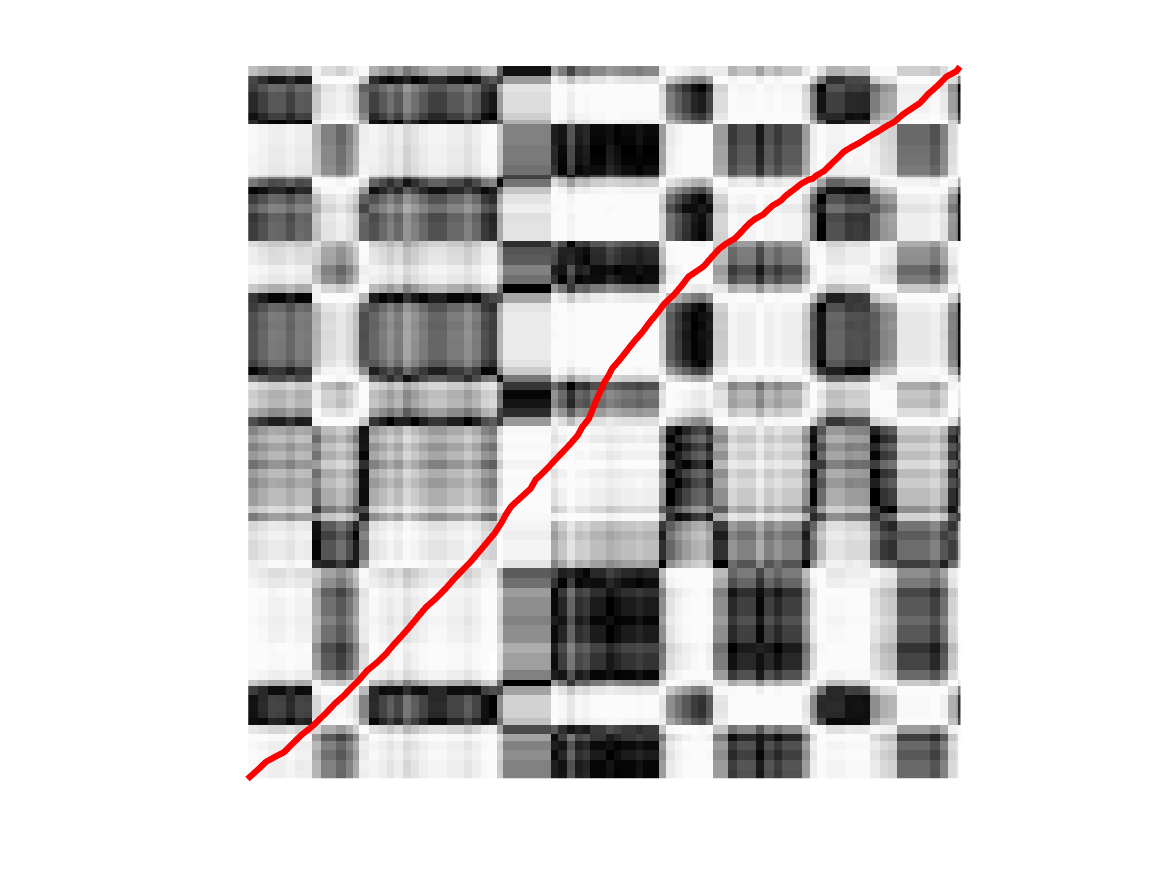}\\
\includegraphics[trim=2cm 0cm 2cm 0cm,clip,width=0.3\textwidth]{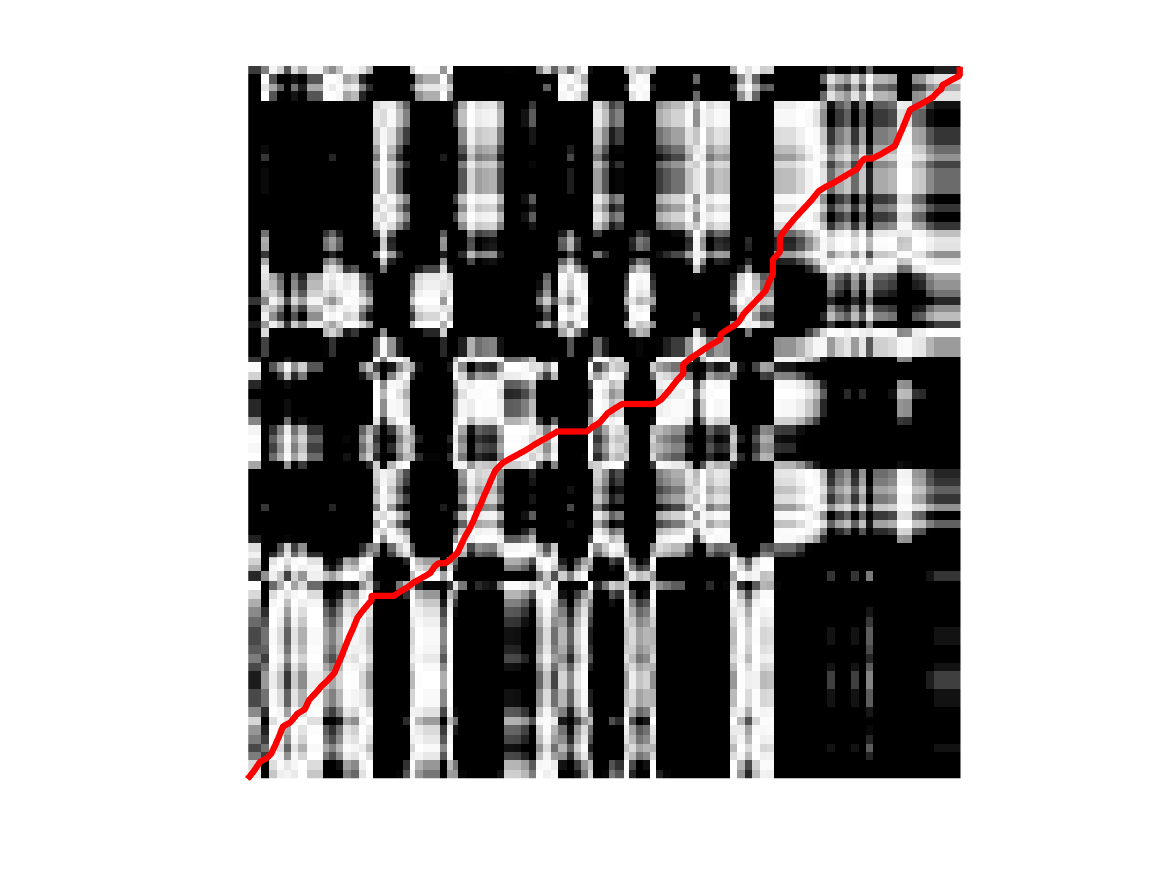}
\includegraphics[trim=2cm 0cm 2cm 0cm,clip,width=0.3\textwidth]{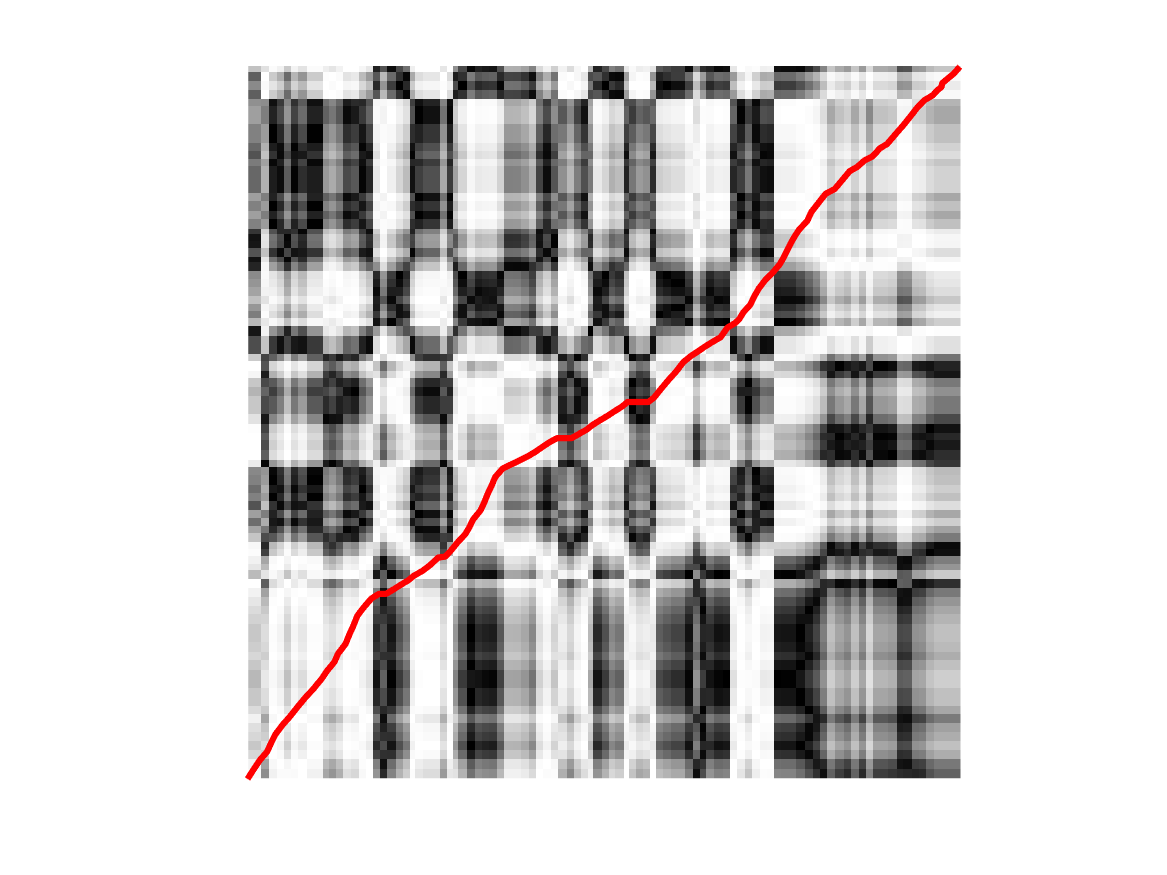}
\includegraphics[trim=2cm 0cm 2cm 0cm,clip,width=0.3\textwidth]{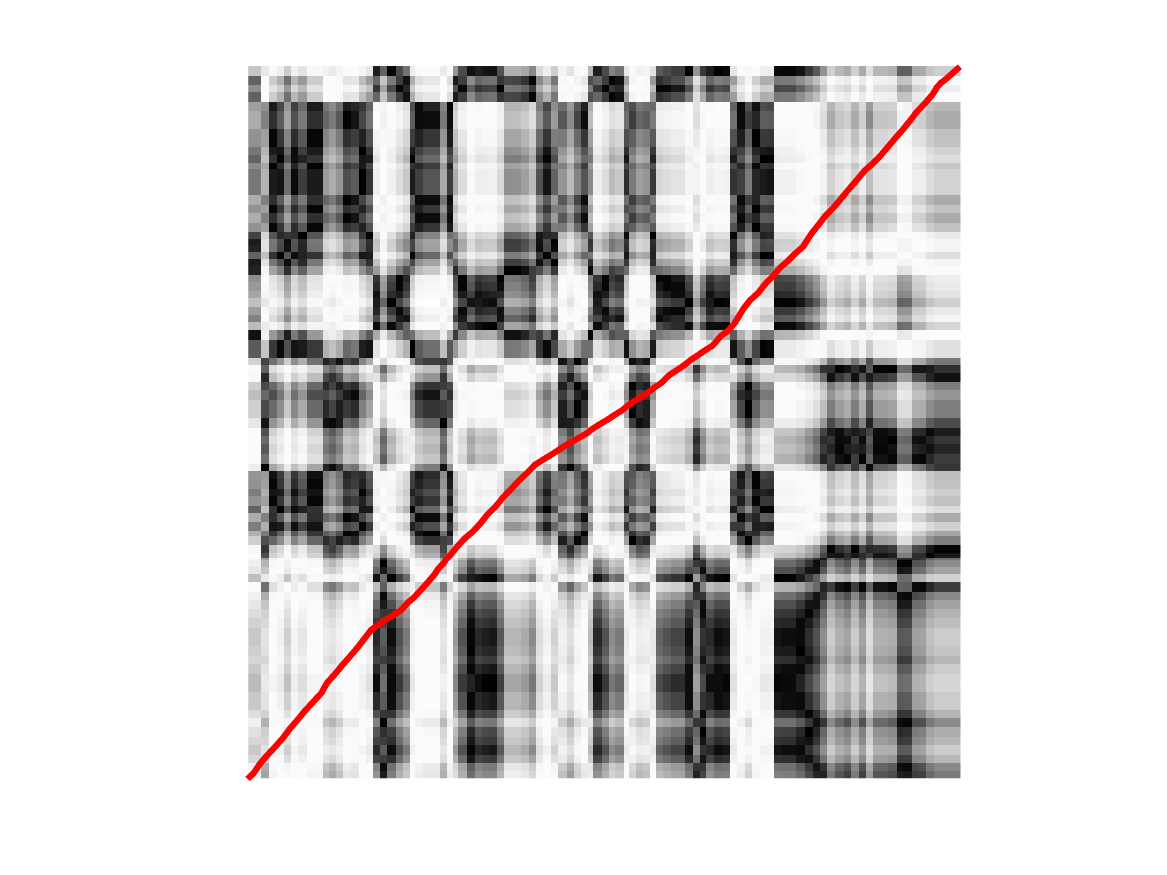}
\caption{Optimal correspondences estimated with the three metamorphoses in Figure \ref{fig:meta.curves}. Background images represent the function $\max(f_\sig,0)$ in \eqref{eq:f.sig},  black color being associated with zeros. The solution attempts to move across higher values of the function, cutting through zeros along vertical or horizontal lines if needed. Each panel corresponds to the same panel in Figure \ref{fig:meta.curves}.
\label{fig:meta.curves.2}}
\end{figure}

\subsection{Case of closed curves}

\label{sec:closed}
The previous developments were obtained assuming that curves were defined over open intervals, and therefore apply mostly to open curves. Closed curves are defined over $T^1$, the open unit interval where the extremities are identified. The boundary condition on $V$, which was $v(0) = v(1) = 0$ for functions defined over $\Om$, now  only requires $v(0) = v(1)$, offering a new degree of freedom, associated with a change of offset, or initial point of the parametrization, represented by the operation $\tau_\de: s \mapsto s+_* \de$ from $T^1$ to itself (where $+_*$ represents addition modulo 1). We restrict our discussion to the two-dimensional case, in which we assume that the compared curves $a_0$ and $a_1$ have angle representations $\th_0$ and $\th_1$.

One can check easily that the distance in Theorem \ref{th:meta.curve.sig} is equivariant through  transformations $\tau_\de$, so that one can define a distance among closed curves that is invariant to rotations and changes of offset by (taking, for example, $\sig = 1$)
\begin{multline}
\label{eq:d.cur.2.inv}
\bar d_{1, \mathrm{rot}}(a_0, a_1) \\
= 2\inf_{\phi, c, \de}  \arccos \int_{T^1} \sqrt{\prt_s \phi(s)} \abs{\cos((\th_0(x+_* \de) - \th_1\circ
\phi(x)-c)/2)} dx.
\end{multline}

Notice that, even when the resulting distance is still attained at a geodesic (or optimal metamorphosis) on the space of functions $a: T^1\to S^1$,  the corresponding curves at intermediate times are not necessarily closed, because the associated closedness condition requires
\[
\int_{T^1} a(s)\, ds = 0
\]
which is not enforced in this approach.
Optimal trajectories are therefore not constrained to consist only of closed curves, and would typically become open for $t\in (0,1)$, even though they start and end with closed curves. The distance in \eqref{eq:d.cur.2.inv} has been applied to obtain the geodesics shown in Figure \ref{fig:meta.curves}, with an extra step  in order to close the intermediate curves for better visualization. This ``closing'' operation simply consisted in replacing $a$ by $\tilde a = (a-\la)/|a-\la|$ where $\la\in\mR^2$ was adjusted so that $\int_\Om \tilde a\, ds = 0$. \\

%
%
%\begin{figure}
%\begin{center}
%\includegraphics[width=5cm]{FIGURES_PDF/curveMatching1.pdf}
%\includegraphics[width=5cm]{FIGURES_PDF/curveMatching2.pdf}\\
%\includegraphics[width=4cm]{FIGURES_PDF/curveMatching3.pdf}
%\includegraphics[width=6cm]{FIGURES_PDF/curveMatching4.pdf}\\
%\end{center}
%\caption[Diffeomorphic matching between curves]{\label{fig:exp.1D} .}
%\end{figure}
%
To correctly define a geodesic distance on spaces of closed curves, one needs to consider the metric induced on the space  $M_c$ of functions $a:T^1 \to S^1$ such that $\int_{T^1} a\,ds = 0$.  This space, however, is not invariant by change of parameters, so that this induced metric is not associated with a metamorphosis (it is the metric induced on the ``submanifold'' of closed curves, for the metamorphosis metric on open curves modulo change of offset and possibly rotations). The expression of this constraint in terms of the $q$ function defined in \eqref{eq:sqvf} is simple in the case $\sig = 1/2$, for which $q = \sqrt{d\phi}\, \al = \sqrt{d \phi}\, a\circ \phi$, so that, after a change of variable,  
\[
 \int_{T^1} a\,ds = \int_{T^1} |q(u)| \, q(u)\, du\,.
 \]
Even in this case, there exists no closed form for the geodesic energy with fixed final reparametrization, but efficient algorithms have been designed to minimize
\[
4\int_0^1 \|\dot q(t, \cdot)\|^2_2\, dt
\]
subject to the constraints that $\|q\|_2^2 = 1$, $\int_{T^1} |q|\,q\,du = 0$, $q(0) = q_0$ and $q(1) = q_1$ (see, for example, \cite{jksj07}). \\

The metric on closed curves also has a nice interpretation in the case $\sig = 1$. In this case, let $f$ and $g$ denote the two coordinates of the representation $q=\mathcal G(\phi, \th)$ multiplied by $\sqrt 2$, i.e., 
$f = \sqrt{2d\phi} \cos\frac{\th}2$ and $g = \sqrt{2d\phi}\sin \frac{\th}{2}$, where $\al = a\circ \phi^{-1} = (\cos\th, \sin\th)$. The closedness constraint, which is
\[
\int_\Om \cos\th\circ\phi^{-1} \,ds = \int_\Om\sin\th\circ\phi^{-1}\,ds = 0,
\]
becomes
\[
\int_\Om d\phi \cos\th\,ds = \int_\Om d\phi \sin\th\,ds = 0
\]
after a change of variables. Writing $\cos\th = \cos^2\frac{\th}{2} - \sin^2\frac{\th}{2}$ and $\sin\th = 2\cos\frac{\th}2\sin\frac{\th}2$, this is equivalent to 
\[
\|f\|_2^2 -\|g\|_2^2 = \scp{f}{g}_2 = 0.
\]
Because $\|f\|_2^2 + \|g\|_2^2 = 2\int_\Om\prt_s\phi = 2$, we find that the constraint is equivalent to $\|f\|_2^2 = \|g\|_2^2 = 1$ and $\scp{f}{g}_2 = 0$, i.e., to $(f, g)$ forming an orthonormal 2-frame in $L^2(\Om)$, and we have
\[
U_1(\phi, \al) = 2 \int_0^1 \left(\|\prt_t f\|^2 + \|\prt_t g\|^2\right)\, ds,
\]
where the left-hand side is two times the geodesic energy of the path $(f, g)$ in the Stiefel manifold $\mathit{St}(\infty, 2)$. Repeating the arguments made in sections \ref{sec:meta.curves.1} or \ref{sec:meta.curves.2} in this setting shows that the optimal metamorphosis with fixed $\phi_1$ is obtained from the shortest length geodesic in $\mathit{St}(\infty, 2)$ connecting the frames $(\sqrt2 \cos\frac{\th^\de_0}2, \sqrt2 \sin\frac{\th^\de_0}2)$ and 
\[
(\ep \sqrt{2\prt_s\phi_1}\cos\frac{\th_1\circ\phi_1}2, \ep \sqrt{2\prt_s\phi_1}\sin\frac{\th_1\circ\phi_1}2)
\]
 where $\th_0$ and $\th_1$ are angle representations of $a_0$ and $a_1$ and the optimization is made over all possible measurable functions $\ep: \Om\to \{-1,1\}$, and all possible offsets $\de$, with the notation $\th_0^\de(s) = \th_0(s +_*\de)$. (The optimization over $\ep$ results from optimizing over all possible angle representations of the two curves.) There is, however, no closed form expression for the geodesic distance on the Stiefel manifold (although equations for geodesics have been described in  \cite{eas98}), and no simple algorithm to solve this optimization problem. Notice that, if one restricts to smooth curves, the search for an optimal $\ep$ is only over constant functions $\ep = \pm 1$ and optimal geodesics can be obtained using a root-finding algorithm over initial conditions of geodesics in $\mathit{St}(\infty,2)$.\\

The rotation-invariant version of the distance also provides an interesting representation, because a rotation acting on curves simply induces  a rotation of the frame $(f, g)$, and the space of such frames modulo rotation now is the Grassmann manifold $\mathit{Gr}(\infty, 2)$ of two-dimensional subspaces of $L^2(\mR)$. The same analysis carries on, the only difference being that one uses now the geodesic distance on the Grassmannian. This geodesic distance can be computed in quasi closed form \cite{ner01}, and is given by $\sqrt{\mathrm{arccos}^2 \la+ \mathrm{arccos}^2\mu}$, where $\la$ and $\mu$ are the singular values of the matrix
\[
\begin{pmatrix}
\scp{f_0}{f_1}_2 & \scp{f_0}{g_1}_2\\
\scp{g_0}{f_1}_2 & \scp{g_0}{g_1}_2
\end{pmatrix}\,,
\]
$(f_0, g_0)$, $(f_1, g_1)$ being orthogonal bases of the two spaces that are compared. This closed form, however, does not lead to a simple version of the distance when one optimizes over changes of sign in $(f_1, g_1)$. More analysis of this framework (in the smooth case), including explicit computations of the geodesic equation and of the scalar curvature can be found in \cite{ymsm07}.

\section{Conclusion}

We have provided in this paper a new view of the first-order Sobolev metric on curves, allowing us to retrieve existing results and obtain new ones. The reduction of an arbitrary $d$-dimensional problem to two dimensions has not, up to our knowledge, been previously proposed in the literature. Neither was the obtention of explicit distances in the non-smooth case, for $\sig \not\in \{1/2, 1\}$. This metric also provided an original  example of metamorphosis, in which the optimal registration is not necessarily diffeomorphic, which led to possible singularities in the optimal solution within a certain range of parameter. This situation is in contrast with existence results for image matching that were obtained, for example, in \cite{ty05}.

\bibliographystyle{plain}
%\nocite{*}
\bibliography{book3}

\begin{thebibliography}{10}

\bibitem{azencott1996distance}
Robert Azencott, Fran{\c{c}}ois Coldefy, and Laurent Younes.
\newblock A distance for elastic matching in object recognition.
\newblock In {\em Pattern Recognition, 1996., Proceedings of the 13th
  International Conference on}, volume~1, pages 687--691. IEEE, 1996.

\bibitem{bauer2017numerical}
Martin Bauer, Martins Bruveris, Philipp Harms, and Jakob M{\o}ller-Andersen.
\newblock A numerical framework for sobolev metrics on the space of curves.
\newblock {\em SIAM Journal on Imaging Sciences}, 10(1):47--73, 2017.

\bibitem{bauer2014constructing}
Martin Bauer, Martins Bruveris, Stephen Marsland, and Peter~W Michor.
\newblock Constructing reparameterization invariant metrics on spaces of plane
  curves.
\newblock {\em Differential Geometry and its Applications}, 34:139--165, 2014.

\bibitem{eas98}
A.~Edelman, T.~A. Arias, and S.~T. Smith.
\newblock The geometry of algorithms with orthogonality contraints.
\newblock {\em SIAM J. Matrix Anal. Appl.}, 20(2):303--353, 1998.

\bibitem{hty09}
D.~R. Holm, A.~Trouv{\'e}, and L.~Younes.
\newblock The euler poincar{\'e} theory of metamorphosis.
\newblock {\em Quarterly of Applied Mathematics}, 2009.
\newblock (to appear).

\bibitem{jksj07}
S.H. Joshi, E.~Klassen, A.~Srivastava, and I.~Jermyn.
\newblock A novel representation for {Riemannian} analysis of elastic curves in
  $r^n$.
\newblock In {\em Proceedings of CVPR'07}, 2007.

\bibitem{ksmj04}
E.~Klassen, A.~Srivastava, W.~Mio, and S.~H. Joshi.
\newblock Analysis of planar shapes using geodesic paths on shape spaces.
\newblock {\em IEEE Trans. Pattern Anal. Mach. Intell.}, 26(3):372--383, 2004.

\bibitem{kurtek2012elastic}
Sebastian Kurtek, Eric Klassen, John~C Gore, Zhaohua Ding, and Anuj Srivastava.
\newblock Elastic geodesic paths in shape space of parameterized surfaces.
\newblock {\em IEEE transactions on pattern analysis and machine intelligence},
  34(9):1717--1730, 2012.

\bibitem{kurtek2018simplifying}
Sebastian Kurtek and Tom Needham.
\newblock Simplifying transforms for general elastic metrics on the space of
  plane curves.
\newblock {\em arXiv preprint arXiv:1803.10894}, 2018.

\bibitem{mm06a}
P.~W. Michor and D.~Mumford.
\newblock Riemannian geometries on spaces of plane curves.
\newblock {\em J. Eur. Math. Soc.}, 8:1--48, 2006.

\bibitem{mm07}
P.~W. Michor and D.~Mumford.
\newblock An overview of the {Riemannian} metrics on spaces of curves using the
  {Hamiltonian} approach.
\newblock {\em Applied and Computational Harmonic Analysis}, 23(1):74--113,
  2007.

\bibitem{my01}
M.~I. Miller and L.~Younes.
\newblock Group action, diffeomorphism and matching: a general framework.
\newblock {\em Int. J. Comp. Vis}, 41:61--84, 2001.
\newblock {Originally published in electronic form in: Proceedings of SCTV 99,
  http://www.cis.ohio-state.edu/~szhu/SCTV99.html}.

\bibitem{mio2007shape}
Washington Mio, Anuj Srivastava, and Shantanu Joshi.
\newblock On shape of plane elastic curves.
\newblock {\em International Journal of Computer Vision}, 73(3):307--324, 2007.

\bibitem{ner01}
Y.~A. Neretin.
\newblock On {Jordan} angles and the triangle inequality in {Grassmann}
  manifolds.
\newblock {\em Geometriae Dedicata}, 86:81--92, 2001.

\bibitem{richardson2013computing}
Casey~L Richardson and Laurent Younes.
\newblock Computing metamorphoses between discrete measures.
\newblock {\em Journal of Geometric Mechanics}, 5(1), 2013.

\bibitem{Richardson2016}
Casey~L. Richardson and Laurent Younes.
\newblock Metamorphosis of images in reproducing kernel hilbert spaces.
\newblock {\em Advances in Computational Mathematics}, 42(3):573--603, Jun
  2016.

\bibitem{samir2012gradient}
Chafik Samir, P-A Absil, Anuj Srivastava, and Eric Klassen.
\newblock A gradient-descent method for curve fitting on riemannian manifolds.
\newblock {\em Foundations of Computational Mathematics}, 12(1):49--73, 2012.

\bibitem{srivastava2016functional}
Anuj Srivastava and Eric~P Klassen.
\newblock {\em Functional and shape data analysis}.
\newblock Springer, 2016.

\bibitem{su2014statistical}
Jingyong Su, Sebastian Kurtek, Eric Klassen, Anuj Srivastava, et~al.
\newblock Statistical analysis of trajectories on riemannian manifolds: bird
  migration, hurricane tracking and video surveillance.
\newblock {\em The Annals of Applied Statistics}, 8(1):530--552, 2014.

\bibitem{ty00}
A.~Trouv{\'e} and L.~Younes.
\newblock Diffeomorphic matching in 1d: designing and minimizing matching
  functionals.
\newblock In D.~Vernon, editor, {\em Proceedings of ECCV 2000}, 2000.

\bibitem{ty01}
A.~Trouv{\'e} and L.~Younes.
\newblock On a class of optimal matching problems in 1 dimension.
\newblock {\em Siam J. Control Opt.}, 39(4):1112--1135, 2001.

\bibitem{ty05}
A.~Trouv{\'e} and L.~Younes.
\newblock Local geometry of deformable templates.
\newblock {\em SIAM J. Math. Anal.}, 37(1):17--59, 2005.

\bibitem{ty05b}
A.~Trouv{\'e} and L.~Younes.
\newblock Metamorphoses through lie group action.
\newblock {\em Found. Comp. Math.}, pages 173--198, 2005.

\bibitem{xie2013parallel}
Qian Xie, Sebastian Kurtek, Huiling Le, and Anuj Srivastava.
\newblock Parallel transport of deformations in shape space of elastic
  surfaces.
\newblock In {\em Proceedings of the IEEE International Conference on Computer
  Vision}, pages 865--872, 2013.

\bibitem{you98}
L.~Younes.
\newblock Computable elastic distances between shapes.
\newblock {\em SIAM J. Appl. Math}, 58(2):565--586, 1998.

\bibitem{you99}
L.~Younes.
\newblock Optimal matching between shapes via elastic deformations.
\newblock {\em Image and Vision Computing}, 17(5):381--389, 1999.

\bibitem{ymsm07}
L.~Younes, P.~Michor, J.~Shah, and D.~Mumford.
\newblock A metric on shape spaces with explicit geodesics.
\newblock {\em Rend. Lincei Mat. Appl.}, 9:25--57, 2008.

\bibitem{you96}
Laurent Younes.
\newblock A distance for elastic matching in object recognition.
\newblock {\em C. R. Acad. Sci. Paris {S{\'e}r.} I Math.}, 322(2):197--202,
  1996.

\end{thebibliography}

\end{document}